\magnification\magstep1
\def\sqr#1#2{{\vcenter{\hrule height.#2pt              %qed
     \hbox{\vrule width.#2pt height#1pt\kern#1pt
     \vrule width.#2pt}
     \hrule height.#2pt}}}
\def\square{\mathchoice\sqr{5.5}4\sqr{5.0}4\sqr{4.8}3\sqr{4.8}3}
\def\qed{\hskip4pt plus1fill\ $\square$\par\medbreak}

%\headline{\hfill May 28, 2009}

\centerline{\bf Degree Complexity of Matrix Inversion}
\bigskip
\centerline{Eric Bedford and Tuyen Trung Truong}
\bigskip
\medskip
\noindent { Abstract.  }
For a $q\times q$ matrix $x=(x_{i,j})$ we let $J(x)=(x_{i,j}^{-1})$ be the Hadamard inverse, which takes the reciprocal of the elements of $x$.  We let $I(x)=(x_{i,j})^{-1}$ denote the matrix inverse, and we define $K=I\circ J$ to be the birational map obtained from the composition of these two involutions.  We consider the iterates $K^n=K\circ\cdots\circ K$ and determine degree complexity of $K$, which is the exponential rate of degree growth $\delta(K)=\lim_{n\to\infty}\left( deg(K^n) \right)^{1/n}$ of the degrees of the iterates.

\bigskip\medskip
\centerline{\S0. Introduction}

\medskip

Let ${\cal M}_q$ denote the space of $q\times q$ matrices with coefficients in ${\bf C}$, and let ${\bf P}({\cal M}_q)$ denote its projectivization.  We consider two involutions on the space of matrices:  $J(x)=(x_{i,j}^{-1})$ takes the reciprocal of each entry of the matrix $x=(x_{i,j})$, and $I(x)=(x_{i,j})^{-1}$ denotes the matrix inverse.  The composition $K=I\circ J$ defines a birational map of ${\bf P}({\cal M}_q)$.

For a rational self-map $f$ of projective space, we may define its $n$th iterate $f^n=f\circ\cdots\circ f$, as well as the degree $deg(f^n)$.  The degree complexity or dynamical degree  is defined as
$$\delta(f):=\lim_{n\to\infty}(deg(f^n))^{1/n}.$$
 In general it is not easy to determine $\delta(f)$, or even to make a good numerical estimate.  Birational maps in dimension 2 were studied in [DF], where a technique was given that, in principle, can be used to determine $\delta(f)$.  This method, however, does not carry over to higher dimension.  In the case of the map  $K_q$,  the dimension of the space and the degree of the map both grow quadratically in $q$, so it is difficult  to write even a small composition $K_q\circ\cdots\circ K_q$ explicitly.   This paper is devoted to determining $\delta(K_q)$.

\proclaim Theorem.  For $q\ge3$, $\delta(K_q)$ is the largest root of the polynomial $\lambda^2-(q^2-4q +2)\lambda+1$.

The map $K$ and the question of determining its dynamical degree have received attention because $K$ may be interpreted as acting on the space of matrices of Boltzmann weights and as such represents a basic symmetry in certain problems of lattice statistical mechanics (see [BHM], [BM]).  In fact there are many $K$-invariant subspaces $T\subset {\bf P}({\cal M}_q)$ (see, for instance, [AMV1] and [PAM]), and it is of interest to know the values of the restrictions $\delta(K|_T)$.  The first invariant subspaces that were considered are ${\cal S}_q$, the space of symmetric matrices, and ${\cal C}_q$, the cyclic (also called circulant) matrices.  The value $\delta(K|_{{\cal C}_q})$ was found in [BV], and another proof of this was given in [BK1].    Angl\`es d'Auriac, Maillard and Viallet  [AMV2] developed numerical approaches to finding $\delta$   and  found approximate values of $\delta(K_q)$ and $\delta(K|_{{\cal S}_q})$ for $q\le 14$.  A comparison of these values with the (known) values of $\delta(K|_{{\cal C}_q})$  led them to conjecture that $\delta(K|_{{\cal C}_q})=\delta(K_q)=\delta(K|_{ {\cal S}_q})$ for all $q$.  

The Theorem above proves the first of these conjectured equalities.  We note that the second equality, $\delta(K|_{{\cal S}_q})=\delta(K_q)$,  involves additional symmetry, which adds another layer of subtlety to the problem.  An example where additional symmetry leads to additional complication has been seen already with the $K$-invariant space ${\cal C}_q\cap {\cal S}_q$: the value of $\delta(K_{{\cal C}_q\cap{\cal S}_q})$ has been determined in [AMV2] (for prime $q$) and [BK2] (for general $q$), and in the general case it depends on $q$ in a rather involved way.
The reason why the cyclic matrices were handled first was that $K|_{{\cal C}_q}$ (see [BV]) and $K|_{{\cal C}_q\cap{\cal S}_q}$ (see [AMV2]) can be converted to maps of the form $L\circ J$ for certain linear  $L$.  In the case of $K|_{{\cal C}_q}$, the associated map is ``elementary'' in the terminology of [BK1], whereas $K|_{{\cal C}_q\cap{\cal S}_q}$ exhibits more complicated singularities, i.e., blow-down/blow-up behavior.   

In contrast, the present paper treats matrices in their general form, so our methods should be applicable to much wider classes of $K$-invariant subspaces.  Our approach is to replace ${\bf P}({\cal M}_q)$ by a  birationally equivalent manifold $\pi:{\cal X}\to{\bf P}({\cal M}_q)$ and consider the induced birational map $K_{\cal X}:=\pi^{-1}\circ K\circ\pi$.  A rational map $K_{\cal X}$ induces a well-defined linear map $K_{\cal X}^*$ on the cohomology group $H^{1,1}({\cal X})$, and the exponential growth rate of degree  is equal to the exponential growth rate of the induced maps on cohomology: 
$$\delta(K)=\lim_{n\to\infty}\left(||(K_{\cal X}^n)^*||_{H^{1,1}({\cal X})} \right)^{1/n}.$$
Our approach is to choose  ${\cal X}$ so that we can determine $(K^n_{\cal X})^*$ sufficiently well.  A difficulty is that frequently $(K^*)^n\ne (K^n)^*$ on $H^{1,1}$.  In the cases we consider, $H^{1,1}$, the cohomology group in (complex) codimension 1, is generated by the cohomology classes corresponding to complex hypersurfaces.  So in order to find a suitable regularization ${\cal X}$, we need to analyze the singularity of the blow-down behavior of $K$, which means that we analyze $K$ at the hypersurfaces $E$ with the property that $K(E)$ has codimension $\ge 2$.

Let us give the plan for this paper.  In general, $deg(K\circ K)\le deg (K)^2$, so $\delta(K)\le deg(K)$.    On the other hand,  $\delta$ decreases when we restrict to a linear subspace, so $\delta(K)\ge \delta(K|_{{\cal C}_q})$.  The paper [BV] shows that $\delta(K|_{{\cal C}_q})$ is the largest root of the polynomial $\lambda^2-(q^2-4q +2)\lambda+1$, so  it will suffice to show that this number is also an upper bound for $\delta(K)$.  In order to find the right upper bound on $\delta(K_q)$, we construct a blowup space $\pi:{\cal Z}\to {\bf P}({\cal M}_q)$.  Such a blowup  induces a birational map $K_{\cal Z}$ of ${\cal Z}$.  Each birational map induces a linear mapping $K_{\cal Z}^*$ on the Picard group $Pic({\cal Z})\cong  H^{1,1}({\cal Z})$.  A basic property is that $\delta(K_{{\cal Z}})\le sp(K_{\cal Z}^*)$, where $sp(K_{\cal Z}^*)$ indicates the spectral radius, or modulus of the largest eigenvalue of $K_{\cal Z}^*$.  Thus the goal of this paper is to construct a space ${\cal Z}$ such that the spectral radius of $K_{\cal Z}^*$ is the number given in the Theorem.  

\bigskip\centerline{\S1.  Basic properties of $I$, $J$, and $K$}
\medskip
For $1\le j\le q-1$, define $R_j$  as the set of matrices in ${\cal M}_q$ of rank less than or equal to $j$.   In ${\bf P}({\cal M}_q)$, $R_1$ consists of matrices of rank exactly 1 since the zero matrix is not in ${\bf P}({\cal M}_q)$.  For $\lambda,\nu\in{\bf P}^{q-1}$, let
$\lambda\otimes\nu=(\lambda_i\nu_j)\in{\bf P}({\cal M}_q)$
denote the outer vector product.  The map
$${\bf P}^{q-1}\times {\bf P}^{q-1}\ni (\lambda,\nu)\mapsto \lambda\otimes\nu \in R_1\subset {\bf P}({\cal M}_q)$$
is biholomorphic, and thus $R_1$ is a smooth submanifold.

We let $I:{\bf P}({\cal M}_q)\to {\bf P}({\cal M}_q)$ denote the birational involution given by matrix inversion $I(A)=A^{-1}$.  We let $x_{[k,m]}$ denote the $(q-1)\times(q-1)$ sub-matrix of $(x_{i,j})$ which is obtained by deleting the $k$-th row and the $m$-th column.  We recall the classic formula $I(x) = (det(x))^{-1}\hat I(x)$, where $\hat I=(\hat I_{i,j})$ is the homogeneous polynomial map of degree $q-1$ given by the cofactor matrix
$$\hat I_{i,j}(x) = C_{j,i}(x)=(-1)^{i+j} det(x_{[j,i]}).\eqno(1.1)$$
Thus $\hat I$ is a homogeneous polynomial map which represents $I$ as a map on projective space.  We see that $\hat I(x)=0$ exactly when the determinants of all $(q-1)\times(q-1)$ minors of $x$ vanish.  

We may always represent a rational map $f=[f_1:\cdots:f_{q^2}]$ of projective space ${\bf P}^{q^2-1}$ in terms of homogeneous polynomials of the same degree and without common factor.  We define the degree of $f$ to be the degree of $f_j$, and the indeterminacy locus is defined as ${\cal I}(f)=\{f_1=\cdots = f_{q^2}=0\}$.  The indeterminacy locus represents the points where it is not possible to extend $f$, even as a continuous mapping.  The indeterminacy locus always has codimension at least 2.  In the case of the rational map $I$,  the polynomials $C_{j,i}(x)$ have no common factor.   Further,  $\hat I(x)=0$ exactly when $x\in R_{q-2}$, so it follows that the indeterminacy set is ${\cal I}(I)=R_{q-2}$.

We let $J:{\bf P}({\cal M}_q)\to {\bf P}({\cal M}_q)$ be the birational involution given by $J(x)=(J(x)_{i,j}) = (1/x_{i,j})$, which takes the reciprocal of all the entries.  In the sequel, we will sometimes write $J(x) = {1\over x}$.  We may define 
$$\hat J(x) =  J(x) \Pi(x)\eqno(1.2)$$ 
where $\Pi(x)=\prod x_{a,b}$ is the homogeneous polynomial of degree $q^2$ obtained by taking the product of all the entries $x_{a,b}$ of $x$, and $\hat J(x)=(\hat J_{i,j})$ is the matrix of homogeneous polynomials of degree $q^2-1$ such that $\hat J_{i,j}=\prod_{(a,b)\ne (i,j)} x_{a,b}$ is the product of all the $x_{a,b}$ except $x_{i,j}$.  Thus $\hat J$ is the projective representation of $J$ in terms of homogeneous polynomials.

We define $K=I\circ J$.  On projective space the map $K$ is represented by the polynomial map (1.4) below.  Since $\hat I\circ\hat J$ has degree $(q-1)(q^2-1)$, we see from Proposition 1.1,  that  the entries of $\hat I\circ \hat J$ must have a common factor of degree $q^3-2q^2$.

 When $V$ is a variety, we write $K(V)=W$ for the strict transform of $V$ under $K$, which is the same as the closure of  $K(V-{\cal I}(K))$.  We say that a hypersurface $V$ is { exceptional} if $K(V)$ has codimension at least 2.  The map $I$ is a biholomorphic map from ${\cal M}_q-R_{q-1}$ to itself, so the only possible exceptional hypersurface for $I$ is $R_{q-1}$.  We define
$$\Sigma_{i,j}=\{ x=(x_{k,\ell}) \in{\cal M}_q: x_{i,j}=0\}.  \eqno(1.3) $$
The map $J$ is a biholomorphic map of ${\cal M}_q-\bigcup_{i,j}\Sigma_{i,j}$ to itself, and the exceptional hypersurfaces are the $\Sigma_{i,j}$.   Further, the indeterminacy locus is
$${\cal I}(J)=\bigcup_{(a,b)\ne(c,d)} \Sigma_{a,b}\cap\Sigma_{c,d}.$$

\proclaim Proposition 1.1.  The degree of $K$ is $q^2-q+1$.  Its representation $\hat K=(\hat K_{i,j})$ in terms of homogeneous polynomials is given by
$$\hat K_{i,j}  (x) = C_{j,i}\left( {1/ x} \right)\Pi(x)\eqno(1.4)$$
where $C_{j,i}$ and $\Pi$ are as in (1.1) and (1.2).

\noindent{\it Proof. }  Observe that $C_{j,i}(1/x)$ is independent of the variable $x_{j,i}$, while $\hat K(x)_{i,j}$ is not divisible by the variables $x_{k,\ell}$ with $k\ne j$ and $\ell\ne i$.  Hence the greatest common divisor of all polynomials on the right hand side of (1.4) is 1.  Thus the algebraic degree of $K$ is equal to the degree of $\hat K(x)_{i,j}$, which is $q^2-q+1$.   \qed

\bigskip
\centerline{\S2.  Construction of ${\cal R}^1$ } 
\medskip
We will construct a complex manifold $\pi:{\cal Z}\to {\bf P}({\cal M}_q)$ by performing a series of blowups.  First we will blow up the spaces $R_1$ and $A_{i,j}$, $1\le i,j\le q$.  The exceptional (blowup) hypersurfaces will be denoted ${\cal R}^1$ and ${\cal A}^{i,j}$. Then we will blow up surfaces $B_{i,j}\subset {\cal A}^{i,j}$, which will create exceptional hypersurfaces ${\cal B}^{i,j}$.  The precise nature of ${\cal Z}$ depends on the order in which the various blowups are performed.  Different orders of blowup will produce different spaces ${\cal Z}$, but the identity map of ${\bf P}({\cal M}_q)$ to itself induces a  birational equivalence between the spaces, and this equivalence induces the identity map on $Pic({\cal Z})$ (as well as on $H^{1,1}({\cal Z})$).  Any of these spaces ${\cal Z}$ yields an induced birational map $K_{\cal Z}$, and each $K_{\cal Z}$ induces essentially the same pullback map $K_{\cal Z}^*$ on $Pic({\cal Z})$.  

We start our discussion with $R_1$.  Let $\pi_1:{\cal Z}_1\to {\bf P}({\cal M}_q)$ denote the blowup of ${\bf P}({\cal M}_q)$ along $R_1$.  We will give a coordinate chart for points of ${\cal Z}_1$ lying over a point  $x^0\in R_1$.  Let us first make a general observation.  Let $\rho_{\ell,m}$ denote the matrix operation which interchanges the $\ell$-th and $m$-th rows of a matrix $x\in{\cal M}_q$, and let $\gamma_{\ell,m}$ denote the interchange of the $\ell$-th and $m$-th columns.  It is evident that $J$ commutes with both $\rho_{\ell,m}$ and $\gamma_{\ell,m}$, whereas we have
$\rho_{\ell,m}(I(x))=I(\gamma_{\ell,m}(x))$.
Thus, for the purposes of looking at the induced map $K_{{\cal Z}_1}$, we may permute the coordinates of $(x_{i,j})$, and without loss of generality we may assume that the (1,1) entry of $x^0$ does not vanish.  This means that we may assume that $x^0=\lambda^0\otimes\nu^0$ with $\lambda^0,\nu^0\in U_1$, where $U_1=\{z=(z_1,\dots,z_q)\in {\bf C}^q:z_1=1\}$.

We write  the standard affine coordinate charts for ${\bf P}({\cal M}_q)$ as
$$W_{r,s}=\{x\in{\cal M}_q: x_{r,s}=1\}\subset {\bf C}^{q^2}, \eqno(2.1)$$
where $1\le r,s\le q$.  Let us define $V$ to be the set of all matrices $x\in{\cal M}_q$ such that the first  row and  column vanish.   Further, for $2\le k,\ell\le q$, we define a subset of $V$:
$$V_{k,\ell}=\{ x\in {\cal M}_q: x = \pmatrix{0 & 0\cr 0 & x_{[1,1]}} {\rm \ and\ } x_{k,\ell}=1 \}.\eqno(2.2)$$
Now we may represent a coordinate neighborhood of ${\cal Z}_1$ over $x^0$ as
$$\pi_1: {\bf C}\times U_1\times U_1\times V_{k,\ell} \to W_{1,1}, \ \ \  \pi_1(s,\lambda,\nu,v) = \lambda\otimes\nu + s v.  \eqno(2.3)$$
Since $\lambda\otimes\nu$ has rank 1 and nonvanishing (1,1) entry, we see that $\pi_1(s,\lambda,\nu,v)\in R_1$ exactly when $s=0$.  Thus the points of ${\cal R}^1$ which are in this coordinate neighborhood are given by $\{s=0\}$.  If $ y\in{\cal M}_q$ is a matrix with $ y_{k,\ell}\ne0$, then we find $\pi_{1}^{-1}(y) = (s,\lambda,\nu,v)$, where
$$\tilde y=y/y_{k,\ell}, \ \ s = y_{k,\ell}, \ \ \lambda=\tilde y_{*,1}, \ \ \nu=\tilde y_{1,*}, \ \ v = s^{-1}(\tilde y-\lambda\otimes\nu).  \eqno(2.4)$$
We may write the induced map $K_{{\cal Z}_1}=\pi_1^{-1}\circ K\circ \pi_1$ in a neighborhood of ${\cal R}^1$ by using the coordinate projections (2.3) and (2.4).  This allows us to show that $K_{{\cal Z}_1}|_{{\cal R}^1}$ has a relatively simple expression:
\proclaim Proposition 2.1.   We have $K_{{\cal Z}_1}({\cal R}^1)=R_{q-1}$, so ${\cal R}^1$ is not exceptional for $K_{{\cal Z}_1}$.  In fact for $z_0=\pi_1(0,\lambda,\nu, v)\in {\cal R}^1$,
$$K_{{\cal Z}_1}(z)=B\pmatrix{0&0\cr 0& I_{q-1}(v')} A \eqno(2.5)$$
where $I_{q-1}$ denotes matrix inversion on ${\cal M}_{q-1}$, and 
$$v'=\left(  {-v_{j,k}\over \lambda_j^2\nu_k^2}\right)_{2\le j,k\le q}, \ A=\pmatrix{ 1 & 0& \cdots &0\cr
-\lambda_2^{-1} & 1 & & \cr
\vdots & & \ddots & \cr
-\lambda_q^{-1} & & & 1\cr},    B = \pmatrix{ 1 & -\nu_2^{-1}& \cdots &-\nu_q^{-1}\cr
0& 1 & & \cr
\vdots & & \ddots & \cr
0 & & & 1\cr}.  \eqno(2.6)$$

\noindent{\it  Proof.  }  Without loss of generality, we work at points $\lambda,\nu\in U_1$ such that $\lambda_j,\nu_k\ne0$ for all $j,k$ and V such that the $v'$ in (2.6) is invertible.  Then 
$$J(\pi_1(s,\lambda,\nu,v)) = {1\over \lambda\otimes\nu} + sv' + O(s^2) = \pi_1(s+O(s^2),\lambda^{-1},\nu^{-1},v'+O(s)). \eqno(2.7)$$
Observe that
$$A\left({1\over\lambda\otimes\nu}\right)B=\pmatrix{1&0\cr0&0}$$
and
$$s\,Av'B=\pmatrix{0&0\cr 0& sA_{[1,1]} v' B_{[1,1]}}.$$ 
Thus
$$\eqalign{
K_{{\cal Z}_1}(z) &= \pi_1^{-1}\circ I\circ J\circ \pi_1(z)\cr
&= \pi_1^{-1} I\left({1\over \lambda\otimes\nu} + sv' + O(s^2)\right)\cr
&= \pi_1^{-1}\left( B\,  I\left( A\left( {1\over \lambda\otimes\nu} + sv' + O(s^2)\right)B\right) A\right)\cr
&=\pi_1^{-1}\left(B\, I\pmatrix{1&0\cr 0&sv'+O(s^2)} A\right),} $$
and the Proposition follows if we let $s\to0$. \qed

Now we will use the identities
$$K_{{\cal Z}_1}\circ J_{{\cal Z}_1}=I_{{\cal Z}_1}, \ \ \ I_{{\cal Z}_1}\circ K_{{\cal Z}_1}=J_{{\cal Z}_1}.$$
\proclaim Proposition 2.2.  We have $K_{{\cal Z}_1}(JR_{q-1})={\cal R}^1$, and thus $JR_{q-1}$ is not exceptional for $K_{{\cal Z}_1}$.

\noindent{\it Proof.  }  For generic $s$, $\lambda$, $\nu$, $v$, and $v'$ as in (2.6), we have (2.7) in the previous Proposition.  Letting $s\to0$, we see that these points are dense in ${\cal R}^1$, and thus $J_{{\cal Z}_1}{\cal R}^1={\cal R}^1.$    Now 
$$\eqalign{K_{{\cal Z}_1}(J(R_{q-1}))  &= I_{{\cal Z}_1}(R_{q-1}) = I_{{\cal Z}_1}(K_{{\cal Z}_1}{\cal R}^1)\cr
& = J_{{\cal Z}_1}({\cal R}^1)={\cal R}^1,\cr}$$
where the second equality in the first line follows from the previous Proposition.
\qed

\bigskip
\centerline{\S3.   Construction of ${\cal A}^{i,j}$}

\medskip

We let $A_{i,j}$ denote the set of $q\times q$ matrices whose $i$-th row and $j$-th columns consist entirely of zeros.  Let $\pi_2:{\cal Z}_2\to{\bf P}({\cal M}_q)$ denote the space obtained by blowing up along all of the the centers $A_{i,j}$ for $1\le i,j\le q$.  As we discussed earlier, it will be immaterial for our purposes what order we do the blowups in.  Let us fix our discussion on $(i,j)=(1,1)$.  The set $A_{1,1}$ is equal to the set $V$ which was introduced in the previous section.  Let us use the notation
$$U=U_{1,r}=\{z\in{\cal M}_q: z = \pmatrix{*&*\cr *&0_{q-1}}, \ z_{1,r}=1\} \eqno(3.1)$$
for the matrices which consist of zeros except for the first row and column, and which are normalized by the entry $z_{1,r}$.  With this notation and with $W_{k,\ell}$, $V_{k,\ell}$ as in (2.1,2), we define the coordinate chart 
$$\pi_2:{\bf C}\times U\times V_{k,\ell}\to W_{k,\ell}\subset {\cal M}_q, \ \ \ \pi_2(s,\zeta,v)=s\zeta + v = \pmatrix {s\zeta & s\zeta\cr s\zeta & v}. \eqno(3.2)$$ 
Coordinate charts of this form give a covering of ${\cal A}^{1,1}$, and $\{s=0\}$ defines the set ${\cal A}^{1,1}$ within each coordinate chart.   If $ x\in{\cal M}_q$, then we normalize to obtain $\tilde x:=x/{x_{k,\ell}}\in W_{k,\ell}$, and
$$\pi_2^{-1}(x)=(s,\zeta,v), \ \ v=\tilde x_{[1,1]},\ s=\tilde x_{1,r}, \ \zeta=(\tilde x-v)/\tilde x_{1,r}. \eqno(3.3)$$ 
We let $K_{{\cal Z}_2}=\pi_2^{-1}\circ K\circ \pi_2$ denote the induced birational map on ${\cal Z}_2$.
\proclaim Proposition 3.1.  For $1\le r,s\le q$, $K_{{\cal Z}_2}(\Sigma_{r,s})={\cal A}^{s,r}$, and in particular $\Sigma_{r,s}$ is not exceptional for $K_{{\cal Z}_2}$.

\noindent{\it Proof. }  As was noted at the beginning of the previous section, it is no loss of generality to assume $(r,s)=(1,1)$ and $2\le k,\ell\le q$.  For generic $x\in{\cal M}_q$, we may use $\hat K$ from (1.4) and define $y$ by
$$\hat K(x) = \Pi(x)\left(C_{j,i}({1\over x}) \right)=y.$$
We write $\pi(\sigma,\zeta,v)=y$, and we next determine $\sigma$, $\zeta$ and $v$.  
Now let us use the notation $s=x_{1,1}$, so $\Pi(x) = s\Pi'(x)$, where $\Pi'$ denotes the product of all $x_{a,b}$ except $(a,b)=(1,1)$.  For $2\le i,j\le q$, we have
$$y_{i,j}=s\Pi'(x)\left ( {1\over s}a_{i,j}(x)+ O(1)\right)$$
with $a_{i,j}(x) = (-1)^{i+j}det((1/x)_{[j,i],[1,1]})$, which gives
$$v_{i,j}=\tilde y_{i,j}=y_{i,j}/y_{k,\ell}=a_{i,j}(x)+O(s), \ \ 2\le i,j\le q.$$
For generic $x$, we may let $s\to0$, and then the value of $v$ approaches $(a_{i,j}(x))/a_{k,\ell}(x)$ which by (1.4) is just $K_{q-1}(x_{[1,1]})$, normalized at the $(k,\ell)$ slot.

The first row and column of $y$ do not involve the (1,1) entry of the matrix $x$, so $y_{1,*}$ and $y_{*,1}$ are divisible by $s$.  By (3.3), we have $\sigma=y_{1,r}/y_{k,\ell}=O(s)$, so we see that $\sigma\to0$ as $s\to 0$.  

An element of the first row of $y$ is given by $y_{1,j}=(-1)^{j+1} det (1/x_{[j,1]})$.  If we expand this determinant into minors along the top row, we have
$$y_{1,j}=\sum_{2\le p\le q} (-1)^{j+1+p} det\left ( (1/x_{[j,1]})_{[1,p]}\right) x_{1,p}^{-1}$$
We use the notation $y_{1,*}$ and $(1/x_{1,*})$ for the vectors $(y_{1,p})_{2\le p\le q}$ and $(1/x_{1,p})_{2\le p\le q}$.  Thus we find
$y_{1,*}=v\,(1/x_{1,*})$.  It is evident that $y_{1,1}=det(1/x_{[1,1]})$.

Now we consider the range of $K$ near ${\cal B}_{1,1}$.  We have seen that $v=K_{q-1}(x_{[1,1]})$, so the values of $v$ are dense in $V_{k,\ell}$.  Now for fixed $v$, we see that the values of $y_{1,*}$ and $y_{*,1}$ span a $2q-2$ dimensional set.  Thus, as we let the values of $x_{1,*}$ and $x_{*,1}$ range over generic values in ${\bf C}^{q-1}\times {\bf C}^{q-1}$, we see that $\zeta$ is dense in $U$.  Thus $K_{{\cal Z}_2}(\Sigma_{1,1}) = {\cal A}^{1,1}$.  
\qed

\bigskip
\centerline{\S4. Construction of ${\cal B}^{i,j}$}
\medskip

For $1\le i,j\le q$, we let $U_{i,j}=\{\zeta\in{\cal M}_q: \zeta_{[i,j]} = 0\}$ to be the set of matrices for which all entries are zero, except on the $i$-th row and $j$-th column.  In the construction of ${\cal A}^{i,j}$, we may consider $U_{i,j}$ (normalized) to be a coordinate chart in the fiber over a point of $A_{i,j}$.  We define the set $B_{i,j}=\{(s,\zeta,v)\in{\cal A}^{i,j}: s=0, \zeta_{i,j}=0\}$, which has codimension 2 in ${\cal Z}_2$, and we let $\pi_3:{\cal Z}_3\to{\cal Z}_2$  be the new manifold obtained by blowing up all the sets $B_{i,j}$.  Let $K_{{\cal Z}_3}$ denote the induced birational map on ${\cal Z}_3$.  As we have seen before, we may focus our attention on the case  $(i,j)=(1,1)$.  Let us use the $(s,\zeta,v)$ coordinate system (3.2) at ${\cal A}^{1,1}$.  Let $U$ be as in (3.1), and set $U'=\{\zeta\in U:\zeta_{1,1}=0\}$.  We define the coordinate projection
$$\pi_3:{\bf C}\times{\bf C}\times U'\times V_{1,1}\to {\bf C}\times U\times V_{1,1}, \ \   \pi(t,\tau,\xi,v) = (s,\zeta,v), \ \ s = t, \zeta = (t\tau,  \xi), v=v,  \eqno(4.1)$$
where the notation $\zeta= (t\tau,\xi)$ means that $\zeta_{1,1}=t\tau$, and $\zeta_{a,b}=\xi_{a,b}$ for all $(a,b)\ne (1,1)$.   Thus ${\cal B}^{1,1}$ is defined by the condition $\{t=0\}$ in this coordinate chart.  Composing the two coordinate projections, ${\cal Z}_3\to {\cal Z}_2$ and ${\cal Z}_2\to {\cal M}_q$, we have
$$\pi: (t,\tau,\xi,v)\mapsto \pmatrix{t^2\tau & t\xi \cr t\xi & v} =x. \eqno(4.2)$$
From (4.2), we see that $\pi^{-1}(x) = (t,\tau,\xi,v)$, where
$$\tilde x = x/x_{\ell,k}, \ v = \tilde x_{[1,1]}, \  t = \tilde x_{1,r}, \ \tau = \tilde x_{1,1}/t^2, \ \xi_{1,j}=x_{1,j}/x_{1,r}, \ 2\le j\le q. \eqno(4.3).$$ 

We will use the following homogeneity property of $K$.  If $x\in{\cal M}_q$, we let $\chi_t(x)$ denote the matrix obtained by multiplying the 1st row by $t$ and then the 1st column by $t$, so the (1,1) entry is multiplied by $t^2$.  It follows that $\chi_t J\chi_t=J$ and $\chi_t\, I\chi_t=I$, so  
$$K\pmatrix {\tau &\xi \cr \xi & v} = \pmatrix{ \tau'&\xi'\cr \xi' & v'}  {\rm \ \ implies\ \ } K\pmatrix {t^2\tau &t\xi \cr t\xi & v} = \pmatrix{ t^{2}\tau'&t^{}\xi'\cr t^{}\xi' & v'}.  \eqno(4.4) $$
\proclaim Proposition 4.1.  For $1\le i,j\le q$, we have $K_{{\cal Z}_3}({\cal B}^{i,j})={\cal B}^{j,i}$, and in particular, ${\cal B}^{i,j}$ is not exceptional.

\noindent {\it Proof. }  As before, we may assume that $(i,j)=(1,1)$.  A point near ${\cal B}^{1,1}$ may be represented in the coordinate chart (4.2) as $\pi(t,\tau,\xi,v) =\pmatrix{t^2\tau & t\xi \cr t\xi & v} =x$.   We define $\tau'$, $\xi'$, and $v'$ by the condition $K\pmatrix{\tau & \xi\cr \xi & v}= \pmatrix{\tau'&\xi'\cr \xi' & v'}$, so $K(x)$ is given by the right hand side of (4.4).  By (4.3),  the coordinates $(t'',\tau'',\xi'',v'')=\pi^{-1}K(x)$ are
$$v''=v/v_{k,\ell}, \ t'' = t \xi'_{1,r}/v'_{k,\ell}, \ \tau''=\tau (v'_{k,\ell}/\xi'_{1,r})^2.$$
From this we see that $t''\to0$ as $t\to0$, which means that $K_{{\cal Z}_3}({\cal B}^{1,1})\subset {\cal B}^{1,1}$.  And since $K$ is dominant on ${\bf P}({\cal M}_q)$, we see that $K_{{\cal Z}_3}({\cal B}^{1,1})$ is dense in ${\cal B}^{1,1}$.  \qed

Next we see how ${\cal A}^{i,j}$ maps under $K_{{\cal Z}_3}$.  A point near ${\cal A}^{1,1}$ may be written in coordinates (3.2) as $(s,\zeta,v)$.  We write $K$ of this point in coordinates (4.1) as $(t,\tau,\xi,w)$. 

\proclaim Proposition 4.2.  For $1\le i,j\le q$, we have  $K_{{\cal Z}_3}({\cal A}^{i,j}) \subset {\cal B}^{j,i} $.  Further, ${dt\over ds}\ne0$ at generic points $(0,\zeta,v)\in{\cal A}^{i,j}$.

\noindent {\it Proof. }  Without loss of generality we assume $(i,j)=(1,1)$.  Let us define $x$ and $y$ as
$$x=\pi_2(s,\zeta,v)=\pmatrix{s\zeta & s\zeta\cr s\zeta & x}, \ \ \ y = \hat K(x) = \Pi(x) C \left({1\over x} \right).$$
For $2\le h,m\le q$ there are polynomials $a_{h,m}(\zeta,v)$ and $b_{h,m}(\zeta,v)$ such that
$$y_{1,1}=s^{2q-1} a_{1,1}(\zeta,v), \ y_{1,m} = s^{2q-2} a_{1,m}(\zeta,v), \ y_{h,m}= s^{2q-3}a_{h,m}(\zeta,v) + s^{2q-2}b_{k,m}(\zeta,v).$$
We have $t=s\, a_{1,r}/a_{k,\ell} + O(s^2)$, so $dt/ds \to a_{1,r}/a_{k,\ell}$ as $s\to 0$.  Thus $dt/ds\ne0$ at generic points of ${\cal A}^{1,1}=\{s=0\}$.  By (4.3), we see that 
$$(t,\tau,\xi,w) \to (0,a_{1,1}a_{k,\ell}/a_{1,r}^2, a_{1,*}/a_{1,r}, a_{[1,1]}/a_{k,\ell})\in{\cal B}^{1,1}$$ 
as $s\to0$.    \qed

\centerline{\S5.  Picard Group $Pic({\cal Z})$}
\medskip

We write ${\cal Z}={\cal Z}_3$ and recall that the Picard group $Pic({\cal Z})$ is the set of divisors modulo linear equivalence.  $Pic({\bf P}({\cal M}_q))=\langle H\rangle$ is generated by any hyperplane $H$.   We will work with the following basis for $Pic({\cal Z})$: 
$$\{H, {\cal R}^1, {\cal A}^{i,j},{\cal B}^{i,j}, 1\le i,j\le q\}. \eqno(5.1)$$ 
Now consider the hypersurface $\Sigma_{i,j}$.  Pulling this back under $\pi_1:{\cal Z}_1\to{\bf P}({\cal M}_q)$, we find 
$$\pi_1^*\Sigma_{i,j} = H_{{\cal Z}_1}=\Sigma_{i,j},$$
where $\Sigma_{i,j}$ on the right hand side denotes the strict transform $\pi^{-1}\Sigma_{i,j}$.  The equality between the strict and total transforms follows because the indeterminacy locus  ${\cal I}(\pi_1^{-1}) = R_1$  is not contained in $\Sigma_{i,j}$.  On the other hand, if we define 
$$T_{i,j}:=\{(a,b):a=i{ \rm\ or\ }b=j\}\eqno(5.2),$$ 
then $\Sigma_{i,j}$ contains $A_{a,b}$ exactly when $(a,b)\in T_{i,j}$.
Thus, pulling back under $\pi_2:{\cal Z}_2\to {\cal Z}_1$, we have
$$\pi^*_2\Sigma_{i,j}=H_{{\cal Z}_2} = \Sigma_{i,j} + \sum_{(a,b)\in T_{i,j}}   {\cal A}^{a,b}.$$
We will next pull this back under $\pi_3:{\cal Z}_3\to{\cal Z}_2$.  For this, we note that $B_{a,b}\subset {\cal A}^{a,b}$, and in addition $B_{i,j}\subset \Sigma_{i,j}$.  Rearranging our answer, we have:
$$\Sigma_{i,j}=H_{\cal Z}- {\cal B}^{i,j} - \sum_{(a,b)\in T_{i,j}} \left( {\cal A}^{a,b} + {\cal B}^{a,b} \right).\eqno(5.3) $$

\proclaim Proposition 5.1.  The class of $JR_{q-1}$ in $Pic({\cal Z})$ is given in the basis (5.1)  by
$$JR_{q-1} = (q^2-q) H - (q-1){\cal R}^1 - (2q-3)\sum_{a,b}{\cal A}^{a,b} - (2q-2)\sum_{a,b}{\cal B}^{a,b}. \eqno(5.4)$$

\noindent{\it Proof. }  The polynomial $P(x):=\Pi(x) det({1\over x})$, analogous to (1.4), is irreducible and has degree $q^2-q$.  Thus $JR_{q-1}=\{P=0\} = (q^2-q)H$ in $Pic({\bf P}({\cal M}_q))$.  Now we pull this back under the coordinate projection $\pi_1$ in (2.3).  That is, we evaluate $P(x)$ for $x=\pi_1(s,\lambda,\nu,v))$.  For $s=0$ and generic $\lambda$, $\nu$, and $v$, the entries of $x=\lambda\otimes\nu + sv$ are nonzero, so $\Pi(x)\ne0$.  We will show $det({1\over x})=\alpha s^{q-1}+\cdots$, where $\alpha\ne0$ for generic $\lambda$, $\nu$, and $v$.  By (2.7), we must evaluate $det(M)$ with 
$M = \lambda^{-1}\otimes\nu^{-1} + sv' + O(s^2)$.  Now we do elementary row and columns such as add $\lambda_j^{-1}\nu$ to the $j$th row, and we do not change the determinant.  In this way, we see that $det(M)$ is equal to $det\pmatrix{1&0\cr 0 & sv'+O(s^2)}= \alpha s^{q-1}+\cdots$.  This means that 
$$(q^2-q)H=\pi^*_1(JR_{q-1})  = JR_{q-1} + (q-1){\cal R}^1\in Pic({\cal Z}_1).$$ 

Now we bring this back to ${\cal Z}_2$ by pulling back under the projection $\pi_2$ defined in (3.2).  In this case, we have $\Pi(\pi_2(s,\zeta,v))= \alpha s^{2q-1} + \cdots$, where $\alpha=\alpha(\zeta,v)\ne0$ for generic $\zeta$ and $v$.  On the other hand, we have $det\pmatrix{s^{-1}\zeta^{-1} &s^{-1}\zeta^{-1} \cr s^{-1}\zeta^{-1}&v^{-1}}=s^{-2}\beta + s^{-1}\gamma+\cdots$, and $\beta(\zeta,v)\ne0$ at generic points.  Thus $P(\pi_2(s,\zeta,v))= c s^{q-3}$, which gives the coefficient $2q-3$ for each ${\cal A}^{i,j}$:
$$(q^2-q)H = JR_{q-1} + (q-1){\cal R}^1 + (2q-3)\sum_{i,j}{\cal A}^{i,j}\in Pic({\cal Z}_2).$$
Pulling back to ${\cal Z}_3$ is similar, except that $\Pi(\pi_3(t,\tau,\xi,v)=\alpha t^{2q}+ \cdots.$  Thus we obtain the coefficient $2q-2$ for ${\cal B}^{i,j}$ in (5.4).
\qed

\bigskip
\centerline{\S6.  The induced map $K_{{\cal Z}}^*$  on $Pic({\cal Z})$ }
\medskip
We define the pullback map on functions by composition  $K^*_{\cal Z}\varphi:=\varphi\circ K_{\cal Z}$.  We may apply $K_{\cal Z}^*$  to local defining functions of a divisor, and  since $K_{\cal Z}$  is well defined off the indeterminacy locus, which has codimension $\ge2$,  $K_{\cal Z}^*$ induces a well-defined pullback  map  on $Pic({\cal Z})$.
\proclaim Proposition 6.1.  $K_{\cal Z}^*$ maps the basis (5.1) according to:
$$\eqalign{ H & \mapsto (q^2-q+1)H -(q-2){\cal R}^1-\sum_{a,b}\left( (2q-3){\cal A}^{a,b}-(2q-2){\cal B}^{a,b}  \right)  \cr
{\cal R}^1 & \mapsto (q^2-q)H-(q-1){\cal R}^1 -  \sum_{a,b}\left( (2q-3){\cal A}^{a,b}-(2q-2){\cal B}^{a,b}  \right)  \cr
{\cal A}^{i,j} & \mapsto H- {\cal B}^{j,i} -  \sum_{(a,b)\in T_{i,j}} \left({\cal A}^{a,b}  +  {\cal B}^{a,b}\right) \cr
{\cal B}^{i,j} & \mapsto {\cal A}^{j,i} + {\cal B}^{j,i} \cr }\eqno(6.1)$$

\noindent{\it Proof. }  Let us start with ${\cal R}^1$.  By \S2, $K_{{\cal Z}}|_{JR_{q-1}}$ is dominant as a map to ${\cal R}^1$.  Since $K_{\cal Z}$ is birational, it is a local diffeomorphism at generic points of $JR_{q-1}$.  Thus we have $K^*_{\cal Z}({\cal R}^1) = JR_{q-1}$, so the second line in (6.1) follows from Proposition 5.1.   

Similarly, since $K_{\cal Z}|_{\Sigma_{i,j}}$ is a dominant map to ${\cal A}^{j,i}$, we have $K_{\cal Z}^*({\cal A}^{i,j})=\Sigma_{j,i}$, and the third line of (6.1) follows from (5.3).

In the case of ${\cal B}^{i,j}$, we know from \S4 that $K_{\cal Z}^{-1}{\cal B}^{i,j}={\cal A}^{j,i}\cup{\cal B}^{j,i}$.  Thus $K_{\cal Z}^*{\cal B}^{i,j} = \lambda{\cal A}^{j,i}+\mu{\cal B}^{j,i}$ for some integer weights $\lambda$ and $\mu$.  Again, since $K_{\cal Z}$ is birational, and $K_{\cal Z}|_{{\cal B}^{i,j}}$ is a dominant map to ${\cal B}^{j,i}$, we have $\mu=1$.   Proposition 4.2 gives us $\lambda=1$.

Finally, set $h(x) = \sum_{i,j} a_{i,j} x_{i,j}$, and let $H=\{h=0\}$ be a hyperplane.  The pullback is given by the class of $\{h\hat K(x)=0\} = \sum_{i,j}a_{i,j}\hat K_{i,j}(x)=0$, where $\hat K$ is given by (1.4).  Pulling back $h$  is similar to the situation in Proposition 5.1, where we pulled back the function $P(x)$.  The difference is that instead of working with $det({1\over x})$ we are working with all of the $(q-1)\times(q-1)$ minors.  By Proposition 1.1, we have $K^*H=(q^2-q+1)H\in Pic({\bf P}({\cal M}_q))$.  Next we will move up to ${\cal Z}_1$ by pulling back under $\pi_1$ and finding the multiplicity of ${\cal R}^1$.  We consider $h\hat K\pi_1(s,\lambda,\nu,v)$, and we recall the matrix $M$ from the proof of Proposition 5.1.  We see that each $(q-1)\times(q-1)$ minor of $M$ is either $O(s^{q-1})$ or $O(s^{q-2})$.  Thus for a generic hyperplane, the order of vanishing is $q-2$, so we have
$$(q^2-q+1)H=K^*H + (q-2){\cal R}^1\in Pic({\cal Z}_1).$$
Next, to move up to ${\cal Z}_2$, we look at the order of vanishing of $h\hat K\pi_2(s,\zeta,v)$ in $s$.  Again $\Pi(\pi_2(s,\zeta,v)) = \alpha s^{2q-1}+\cdots$.  The $(q-1)\times(q-1)$ minors of 
$\pmatrix{s^{-1}\zeta^{-1} &s^{-1}\zeta^{-1} \cr s^{-1}\zeta^{-1}&v^{-1}}$  
which grow most quickly behave like $s^{-2}\beta + s^{-1}\gamma+\cdots$.  Thus for generic coefficients $a_{i,j}$ we have vanishing to order $2q-3$ in $s$, and so $2q-3$ is the coefficient for each ${\cal A}^{i,j}$ as we pull back to $Pic({\cal Z}_2)$.  Coming up to ${\cal Z}_3={\cal Z}$, we pull back under $\pi_3$, and the calculation of the multiplicity of ${\cal B}^{i,j}$ is similar.  This gives  the first line in (6.1).  \qed

\proclaim Proposition 6.2.  The characteristic polynomial of the transformation (6.1) is 
$$P(\lambda)Q(\lambda)^{q-1}(\lambda-1)^{q^2-q+2}(\lambda+1)^{q^2-3q+2},$$ 
where
$P(\lambda)=\lambda^2-(q^2-4q+2)\lambda+1$ and $Q = (\lambda^2+1)^2-(q-2)^2\lambda^2$.

\noindent{\it Proof. }  We will exhibit the invariant subspaces of $Pic({\cal Z})$ which correspond to the various factors of the characteristic polynomial.  First, we set ${\cal A}:=\sum{\cal A}^{k,\ell}$ and ${\cal B}:=\sum{\cal B}^{k,\ell}$, where we sum over all $k$ and $\ell$, and we set $S_1=\langle H,{\cal R}^1,{\cal A},{\cal B}\rangle$.  By (6.1), $S_1$ is $K_{{\cal Z}}^*$-invariant, and the characteristic polynomial of $K_{\cal Z}^*|_{S_1}$ is seen to be $P(\lambda)(\lambda-1)^2$. 

Next, if $i<j$, then we set $\alpha_{i,j}={\cal A}^{i,i}+{\cal A}^{j,j}-({\cal A}^{i,j}+{\cal A}^{j,i})$, and similarly for $\beta_{i,j}$, using the ${\cal B}^{k,\ell}$.  Then by (6.1), $S_{i,j}:=\langle\alpha_{i,j},\beta_{i,j}\rangle$ is invariant, and the characteristic polynomial of $K_{\cal Z}^*|_{S_{i,j}}$ is $(\lambda-1)^2$.

Similarly, if $i<j<k$, we set $\alpha_{i,j,k}={\cal A}^{i,i}+{\cal A}^{j,j}+{\cal A}^{k,k}-\left( {\cal A}^{i,j}+{\cal A}^{j,k} + {\cal A}^{k,i}\right)$ and define $\beta_{i,j,k}$ similarly.  Then the 2-dimensional subspace  $S_{i,j,k}:=\langle \alpha_{i,j,k},\beta_{i,j,k}\rangle$ is invariant, and the characteristic polynomial of  $K_{\cal Z}^*|_{S_{i,j,k}}$ is  $(\lambda+1)^2$.

Finally, for each $i$, we consider the row and column sums ${\cal A}_{r_i}=q\sum_{j} {\cal A}^{i,j}-{\cal A}$, ${\cal A}_{c_j}=q\sum_{i} {\cal A}^{i,j}-{\cal A}$, and we make the analogous definition for ${\cal B}_{r_i}$ and ${\cal B}_{c_j}$.  The 4-dimensional subspace $\langle {\cal A}_{r_i}, {\cal A}_{c_i},{\cal B}_{r_i}, {\cal B}_{c_i}\rangle$ is invariant and  yields the factor $Q(\lambda)$.  These invariant subspaces span $Pic({\cal Z})$, and the product of these factors gives the characteristic polynomial stated above. \qed

\noindent{\it Proof of the Theorem. }  The spectral radius of $K^*_{\cal Z}$ is the largest root of the characteristic polynomial, which is given in Proposition 6.2.  By inspection, the largest root of the characteristic polynomial is the largest root of $P(\lambda)$.  The spectral radius of $K^*_{\cal Z}$ is an upper bound for $\delta(K)$.  On the other hand, it was shown in [BV] that this same number is also a lower bound for $\delta(K)$, so the Theorem is proved.  \qed

\bigskip
\centerline{References } 
\medskip
%\item{[AABHM]}  N. Abarenkova, J.-C. Angl\`es d'Auriac, S. Boukraa, S.
%Hassani, and J.-M. Maillard, Rational dynamical zeta functions for
%birational transformations, Phys.\ A 264 (1999), 264--293.

%\item{[AABM]} N. Abarenkova, J.-C. Angl\`es d'Auriac, S. Boukraa,  and
%J.-M. Maillard,  Growth-complexity spectrum of some discrete dynamical
%systems, Physica D 130  (1999), 27--42.

\item{[AMV1]} J.C. Angl\`es d'Auriac, J.M. Maillard, and C.M. Viallet, A classification of four-state spin edge Potts models, J.\ Phys.\ A 35 (2002), 9251--9272.  cond-mat/0209557
\item{[AMV2]}  J.C. Angl\`es d'Auriac, J.M. Maillard, and C.M. Viallet,  On the complexity of some birational transformations. J. Phys. A 39 (2006), no. 14, 3641--3654.  math-ph/0503074
\item{[BK1]}  E. Bedford and K-H Kim, On the degree growth of birational mappings in higher dimension, J. Geom.\ Anal.\ 14 (2004), 567--596.  arXiv:math.DS/0406621
\item{[BK2]}  E. Bedford and K-H Kim, Degree growth of matrix inversion: birational maps of symmetric, cyclic matrices. Discrete Contin. Dyn. Syst. 21 (2008), no. 4, 977--1013. 
%\item{[BMV]} M.P. Bellon, J.-M. Maillard, and C.-M. Viallet, Integrable Coxeter groups, Phys. Lett. A 159 (1991), 221--232.
\item{[BV]} M. Bellon and C.M. Viallet, Algebraic entropy, Comm.\ Math.\ Phys., 204 (1999), 425--437.
%\item{[BTR]} M. Bernardo, T.T. Truong and G. Rollet, The discrete
%Painlev\'e I equations: transcendental integrability and asymptotic
%solutions,  J. Phys. A: Math. Gen., 34 (2001), 3215--3252.
\item{[BHM]} S. Boukraa, S. Hassani, J.-M. Maillard,  Noetherian mappings,
Physica D, 185 (2003), no. 1, 3--44. 
\item{[BM]}  S. Boukraa and J.-M. Maillard, Factorization properties of
birational mappings, Physica A 220 (1995), 403--470.
\item{[DF]}  J. Diller and C. Favre,  Dynamics of birational maps of surfaces, Amer.\ J. Math. 123 (2001), no. 6, 1135--1169.
\item{[PAM]}  E. Preissmann J.-Ch.\ Angl\`es d'Auriac and J.-M. Maillard, Birational mappings and matrix subalgebra from the chiral Potts model, J. of Math.\ Physics 50, 013302 (2009).  arXiv: 0802.1329
%\item{[T]}  T. Truong,  Degree complexity of matrix inversion: symmetric case, in preparation.

\bigskip
\rightline{\tt bedford@indiana.edu}

\rightline{\tt truongt@umail.iu.edu}

\rightline{Department of Mathematics}

\rightline{Indiana University}

\rightline{Bloomington, IN 47405}

\bye